\documentclass[12pt,reqno]{amsart}

\numberwithin{equation}{section} 

\usepackage{enumerate}
\usepackage{amssymb}
\usepackage{amsmath}
\usepackage{amscd}
\usepackage{amsthm}
\usepackage{amsfonts}
\usepackage{graphicx}
\usepackage[all]{xy}
\usepackage{verbatim}
\usepackage{hyperref}

\newtheorem{theorem}{Theorem}[section]

\newtheorem{lemma}[theorem]{Lemma}
\newtheorem{corollary}[theorem]{Corollary}
\newtheorem{definition}[theorem]{Definition}

\theoremstyle{definition}

\theoremstyle{definition}
\newtheorem{question}[theorem]{Question}
\theoremstyle{definition}\newtheorem{remarks}[theorem]{Remarks}
\theoremstyle{definition}\newtheorem{remark}[theorem]{Remark}
\theoremstyle{definition}\newtheorem*{acknowledgments}{Acknowledgments}

\newcommand{\al}{\alpha}
\newcommand{\be}{\beta}
\newcommand{\ga}{\gamma}
\newcommand{\Ga}{\Gamma}
\newcommand{\del}{\delta}

\newcommand{\lam}{\lambda}
\newcommand{\Lam}{\Lambda}
\newcommand{\eps}{\epsilon}

\newcommand{\ka}{\kappa}

\newcommand{\vphi}{\varphi}

\newcommand{\cM}{\mathcal{M}}

\newcommand{\bA}{\mathbb{A}}

\newcommand{\bC}{\mathbb{C}}

\newcommand{\bG}{\mathbb{G}}

\newcommand{\bR}{\mathbb{R}}
\newcommand{\bZ}{\mathbb{Z}}
\newcommand{\bQ}{\mathbb{Q}}

\newcommand{\bN}{\mathbb{N}}
\newcommand{\bH}{\mathbb{H}}
\newcommand{\bT}{\mathbb{T}}

\newcommand{\SL}{\operatorname{SL}}
\newcommand{\PGL}{\operatorname{PGL}}
\newcommand{\GL}{\operatorname{GL}}
\newcommand{\PSL}{\operatorname{PSL}}

\newcommand\norm[1]{\left\|#1\right\|}
\newcommand\set[1]{\left\{#1\right\}}
\newcommand\pa[1]{\left(#1\right)}

\newcommand\idist[1]{\langle#1\rangle}
\newcommand\av[1]{\left|#1\right|}

\newcommand\diag[1]{\operatorname{diag}\left(#1\right)}

\newcommand{\onto}{\xymatrix{\ar@{>>}[r]&}}
\newcommand{\da}[4]{\xymatrix{#1 \ar@<.5ex>[r]^{#2} \ar@<-.5ex>[r]_{#3} & #4}}

\begin{document}
\title{Diophantine approximations on fractals}
\author{M. Einsiedler, L. Fishman, U. Shapira.}
\thanks{Part of the third author's Ph.D thesis at the Hebrew University of Jerusalem. The partial support of ISF grant number 1157/08, as well as the  Binational Science Foundation
grant 2008454 is acknowledged.}
\begin{abstract}
We exploit dynamical properties of diagonal actions to derive results in Diophantine approximations. In particular, we prove that the continued fraction expansion of almost any point on the middle third Cantor set (with respect to the natural measure)  contains all finite patterns (hence is well approximable). Similarly, we show that for a variety of fractals in $[0,1]^2$, possessing some 
symmetry, almost any point is not Dirichlet improvable (hence is well approximable) and has property C (after Cassels).
We then settle by similar methods a conjecture of M. Boshernitzan saying that there are no irrational numbers $x$ in the unit interval such that the continued fraction expansions of $\set{nx \textrm{ mod } 1}_{n\in\bN}$ are uniformly eventually bounded. 
\end{abstract}
\maketitle
\section{Introduction}\label{introduction}
\subsection{Preface}\label{Preface}
In  the theory of metric Diophantine approximations, one wishes to understand how well vectors in $\bR^d$ can be
approximated by rational vectors. The quality of approximation can be measured in various forms leading to numerous
Diophantine classes of vectors such as \textit{WA (well approximable), VWA (very well approximable), DI (Dirichlet improvable)} and so forth. Usually such a class is either a null set or generic (i.e.\ its complement is a null set) and often one encounters the phenomena of the class being null but of full dimension.
Given a closed subset $M\subset\bR^d$ supporting a natural measure (for example a lower dimensional submanifold with the volume measure or a fractal with the Hausdorff measure),
it is natural to investigate the intersection of $M$ with the various Diophantine classes.
It is natural to  expect that unless there are  obvious obstacles, the various Diophantine classes will intersect $M$ in a set which will inherit the characteristics of the class, i.e.\ if the class is null, generic or of full dimension in $\bR^d$, then its intersection with $M$ would be generic, null or of full dimension in $M$ as well. 

Let us demonstrate this with two examples in the real line. We consider the intersection of the middle third
Cantor set, $C$, in the unit interval, with two classes: BA and VWA. The class BA of badly approximable numbers  consists of  real numbers whose coefficients in their continued fraction
expansion are bounded and the class WA is its complement. A classical result of Borel says that BA is null. Nevertheless, Schmidt showed in ~\cite{Sch} that 
it is of dimension 1. 
It was shown independently  in ~\cite{KWIJM} and ~\cite{KTV}, that the dimension of $C\cap \textrm{BA}$ is full, i.e.\ equals
$\log 2/\log 3$ (see ~\cite{F} for some recent developments). One of the motivating 
questions for this paper, answered affirmatively in Corollary ~\ref{cantor set and WA}, was to decide whether
$C\cap \textrm{BA}$ is null with respect to the Hausdorff measure on $C$. \\
The class VWA (in the real line) is a subclass of WA and consists of numbers $x$ for which there exists $\del >0$ such that one can find infinitely many solutions over $\bZ$ to the inequality $\av{qx-p}<q^{-(1+\del)}$. It is null and of full dimension in $\bR$. It was shown  in ~\cite{W1} that $C\cap\textrm{VWA}$ is null with respect to the Hausdorff
Measure on $C$ and in ~\cite{LSV} a lower bound for the dimension of this intersection is given. As far as we know 
 it is not known if the dimension equals $\dim C$.

The intersection of the class of VWA vectors with submanifolds and fractals in $\bR^d$ has attracted much 
 attention. As this class will not concern 
us in this paper we refer the reader to the following papers for further discussions: ~\cite{KLW},  
~\cite{KM},~\cite{W1}, \cite{K},\cite{PV}, and  \cite{LSV}.

In this paper we will be concerned with inheritance of genericity to certain fractals in $\bR$ and $\bR^2$,
with respect to three Diophantine classes WA, DI and C (see Definitions 
~\ref{def WA}, ~\ref{def improvable},~\ref{property C}).
In Theorems ~\ref{theorem 1} through \ref{theorem 2'} we prove that the above Diophantine classes remain generic or null
when additional assumptions on the fractal (and the measure supported on it) are imposed. These involve positivity
of dimension and invariance under an appropriate map.
The reader is referred to Remarks ~\ref{remark 1} for a discussion about the necessity of these
additional assumptions as well as the restriction to dimensions 1 and 2. 

Our arguments rely on the measure classification results obtained by E. Lindenstrauss in ~\cite{L} and by 
M. Einsiedler, E. Lindenstrauss and A. Katok  in ~\cite{EKL}.
\subsection{Diophantine classes}\label{WA}
Vectors in $\bR^d$ will be thought of as column vectors and the action of matrices on them will be from the left. We now define the Diophantine classes we will consider.
\begin{definition}\label{def WA}
A vector $v\in\bR^d$ is said to be well approximable (WA), if for any $\eps>0$ one can find $\vec{m}\in\bZ^d,n\in\bN$ such that 
\begin{equation}\label{WA'}
\norm{nv-\vec{m}}_\infty<\frac{\eps}{n^{1/d}}.
\end{equation} 
We denote $\textrm{WA}=\set{v\in\bR^d:v\textrm{ is WA}}$.
\end{definition}
It is well known that WA is a generic class.
\begin{definition}\label{def improvable}
A vector $v\in\bR^d$ is said to be Dirichlet improvable if there exists $0<\mu<1$, such that for all 
sufficiently large $N$ the following statement holds:
$$\textrm{There exists } \vec{m}\in\bZ^d, n\in\bN\textrm{ such that } 0<n\le  N^d \textrm{ and }
\norm{nv - \vec{m}}<\mu N^{-1}.$$
We denote $\textrm{DI}=\set{v\in\bR^2:v\textrm{ is Dirichlet improvable}}$. We say 
that $v$ is not DI if $v\notin DI$.
\end{definition}
In ~\cite{DS1} Davenport and Schmidt introduced the notion of Dirichlet improvable vectors and showed amongst other things that the class, BA, of badly approximable vectors (which is the complement of WA) is contained in the class DI. Moreover, they showed that in dimension 1 the two classes are equal (modulo the rationals). In ~\cite{DS} it is shown that DI is a null set. Recently
N. Shah, motivated by the work of Kleinbock and Weiss ~\cite{KW}, 
showed in ~\cite{Shah} that the intersection of DI with any non degenerate analytic curve in $\bR^d$ is null as well. 
In the following definition we use the notation $\idist{\ga}$ for the distance of a real number $\ga$ to the integers. 
\begin{definition}\label{property C}
A vector $v\in\bR^d$ is said to have property C (after Cassels) of the first type if the following statement holds:
$$\textrm{For all }\vec{\ga}\in\bR^d \quad \liminf_{\av{n}\to\infty}\av{n}\prod_1^d\idist{nv_i-\ga_i}=0.$$
It is said to have property C of the second type if the following statement holds: 
$$\textrm{For all }\ga\in\bR \quad \liminf_{\vec{n}\in\bZ^d, \prod\av{n_i}\to\infty}\pa{\prod_1^d\av{n_i}}\idist{\sum_1^d n_iv_i -\ga}=0.$$
We denote $\textrm{C}=\set{v\in\bR^d:v\textrm{ has property C of the first and the second type}}$. We say that 
$v$ has property C if $v\in C$.
\end{definition}
For $d=1$, it is shown in \cite{D} that there are no real numbers with property C (of the first or second type which coincide in this case). In ~\cite{Shap1} the third named author showed that for $d\ge 2$ the class C is generic. 
We remark that if a vector $(\al,\be)^t\in\bR^2$ ($t$ stands for transpose), has property C of the first type, 
then in particular, $\al,\be$ satisfy the well known Littlewood conjecture, i.e.\ $$\liminf_{n\to\infty}n\idist{n\al}\idist{n\be}=0.$$
Note that WA, DI and C, are invariant under translations by integer vectors, hence define subsets of the $d$-torus $\bR^d/\bZ^d$. We use the same notation for the corresponding subsets of the $d$-torus.
\subsection{Statements of results}\label{statements of results}
Before stating the main results which we prove in this paper we need to recall the notion of 
\textit{dimension of a measure}. Let $K$ be  a compact metric space and $\mu$ a Borel probability measure 
on $K$. The upper and lower local dimension functions of $\mu$ are defined to be 
\begin{equation}\label{ldf}
\underline{d}_\mu(x)  =\liminf_{r\to 0}\frac{\log\pa{\mu\pa{B_r(x)}}}{\log\pa{r}},\;
\overline{d}_\mu(x)=\limsup_{r\to 0}\frac{\log\pa{\mu\pa{B_r(x)}}}{\log\pa{r}}.
\end{equation} 
We say that $\mu$ has \textit{exact dimension} if there exists a number $d$ such that $\underline{d}_\mu(x)=\overline{d}_\mu(x)=d$ for $\mu$-almost any $x$. In this case we sometime simply say
that $\mu$ is of dimension $d$. 
\begin{remark}\label{dimension is entropy}
In the following theorems the fact that measures have exact dimension follows from the other assumptions;
it follows from~\cite{BK} that both the lower and upper dimension functions are equal almost surely to a positive multiple of the entropy of the system.
\end{remark}
\begin{theorem}\label{theorem 1}
Let $n\in\bN$ and let $\mu$ be a probability measure on the unit interval which is invariant and ergodic under $\times n$ modulo $1$ (i.e.\ under multiplication by $n$ modulo 1), and has positive dimension. Then $\mu$ almost any $x\in [0,1]$ is WA.
\end{theorem}
\begin{theorem}\label{theorem 2}
Let $\ga:\bR^2/\bZ^2\to\bR^2/\bZ^2$ be a hyperbolic automorphism, induced by the linear action of a matrix $\ga\in \SL_2(\bZ)$ and let $\mu$ be a probability measure which is invariant and ergodic with respect to $\ga$, and has positive dimension. Then $\mu$ almost any $v\in \bR^2/\bZ^2$ is WA.
\end{theorem} 
One way to construct examples of measures $\mu$ on the unit interval or on the 2-torus satisfying the assumptions of Theorems ~\ref{theorem 1} and ~\ref{theorem 2} respectively is to choose an appropriate partition of the underlying space for which the resulting factor map to the symbolic system is an isomorphism of measurable dynamical systems. In the case of the unit interval, one chooses the 
partition into $n$ intervals of equal length, and in the case of the 2-torus, a Markov partition corresponding to
 $\ga$ (see \cite{AW}). Then, one takes a (topologically transitive) subshift of finite type of the symbolic system and the unique maximal entropy probability measure 
supported on it and translates this measure to the original space. 

Before stating further results we briefly introduce some notation (see \S~\ref{space of lattices} for a more thorough account). 
For any positive integer $d$ ($d=1$ or $2$ in our discussions), let $\set{a_t}_{t\in\bR},\set{u_v}_{v\in\bR^d}<\PGL_d(\bR)$ be the subgroups given by
\begin{equation}\label{a_t and u_v}
a_t=\operatorname{diag}(e^t,\dots,e^t,e^{-dt}),\quad 
u_v=
\pa{
\begin{array}{ll}
I_d& -v\\
0& 1
\end{array}
},
\end{equation}
where $I_d$ is the $d\times d$ identity matrix. Our arguments rely on the natural identification of the 
$d$-torus with the periodic orbit $\set{u_v\PGL_{d+1}(\bZ) : v\in\bR^d}$ in the homogeneous space $\PGL_{d+1}(\bR)/\PGL_{d+1}(\bZ)$ (see~\eqref{d-torus as periodic orbit}). This enables us to view measures supported on the $d$-torus as measures supported in the space $\PGL_{d+1}(\bR)/\PGL_{d+1}(\bZ)$.

In the following two theorems we are able to obtain stronger results than in Theorems ~\ref{theorem 1},~\ref{theorem 2}. The price is reflected in the stronger assumptions which are automatically satisfied in many applications (see Remark ~\ref{friendly}).
\begin{theorem}\label{theorem 1'}
Let $n$ and $\mu$ be as in theorem ~\ref{theorem 1}. Viewing $\mu$ as a probability measure on 
$\PGL_{2}(\bR)/\PGL_{2}(\bZ)$, if we further assume
that any weak$^*$ limit of $\frac{1}{T}\int_0^T(a_t)_*\mu dt$ is a probability  measure on $ \PGL_2(\bR)/\PGL_2(\bZ)$ (i.e.\ there is no escape of mass on average), then for $\mu$ almost any $s\in \bR/\bZ$ 
\begin{equation}\label{density1}
\overline{\set{a_tu_s\PGL_2(\bZ)}}_{t\ge 0}=\PGL_2(\bR)/\PGL_2(\bZ).
\end{equation} 
Furthermore, if for a given $s$ \eqref{density1} holds, then the continued fraction expansion of $s$ contains all patterns.
\end{theorem} 
\begin{theorem}\label{theorem 2'}
Let $\ga,\mu$ be as in Theorem ~\ref{theorem 2}. Viewing $\mu$ as a probability measure on 
$\PGL_{3}(\bR)/\PGL_{3}(\bZ)$, if we further assume
that any weak$^*$ limit of $\frac{1}{T}\int_0^T(a_t)_*\mu dt$ is a probability  measure on $ \PGL_3(\bR)/\PGL_3(\bZ)$ (i.e.\ there is no escape of mass on average), then for $\mu$ almost any $v\in \bR^2/\bZ^2$ 
$$\overline{\set{a_tu_v\PGL_3(\bZ)}}_{t\ge 0}=\PGL_3(\bR)/\PGL_3(\bZ).$$ 
In particular $v$ is not DI (hence is WA) and has property C.
\end{theorem}
\begin{remark}\label{friendly} In ~\cite{KLW}, Kleinbock, Lindenstrauss, and Weiss showed that if $\mu$ is \textit{friendly} (see \S2 of ~\cite{KLW} for the definition), then there is no escape of mass on average and so the further assumptions in Theorems ~\ref{theorem 1'}, ~\ref{theorem 2'} are satisfied automatically. For a detailed proof of this statement the reader is further referred to ~\cite{R} Corollary 3.2.
\end{remark}

In \S2 of the paper ~\cite{KLW} it is shown that if $F\subset [0,1]$ is a fractal defined as the attractor of an irreducible system of contracting self similar maps satisfying the open set condition, then the Hausdorff measure
on $F$ is of positive dimension and is friendly. In many examples, the fractal $F$ is invariant and ergodic (with respect to the Hausdorff measure) under $\times n$ (for a suitable choice of $n$), hence Theorem ~\ref{theorem 1'} applies by the above remark. In particular we have the following corollary
which served as one of the motivating questions for this work. 
\begin{corollary}\label{cantor set and WA}
Almost any point in the middle third Cantor set (with respect to the natural measure) is WA and moreover its continued fraction expansion contains all patterns.
\end{corollary} 

Our last theorem is of a different nature as it is an everywhere statement. It was conjectured to hold by M. Boshernitzan and communicated by the second named author. 
\begin{theorem}\label{theorem 3}
If we denote for $x\in [0,1]$, $c(x)=\limsup a_n(x)$ where $a_n(x)$ are the coefficients in the continued  fraction expansion of $x$, then for any irrational $x\in [0,1]$, $\sup_nc(nx)=\infty$, where $nx$ is calculated modulo 1.
\end{theorem}
The proofs of all the above theorems are similar in nature. We shall first prove Theorems ~\ref{theorem 2},~\ref{theorem 2'} which are somewhat simpler but contain the ideas. We then prove Theorems ~\ref{theorem 1},~\ref{theorem 1'} which involves $S$-arithmetic arguments and finally prove Theorem ~\ref{theorem 3} which involves adelic arguments. 
\begin{remarks}\label{remark 1}\quad
\begin{enumerate}
\item We note that there are fractals of positive dimension which intersect the various generic Diophantine 
classes trivially. For example it is not hard to construct a closed set of positive dimension in the unit interval which is contained in the class
BA of badly approximable numbers. Hence in order to obtain results as above one must impose some further assumptions, which in our case, are reflected in the 
symmetry of the fractal given by the invariance under the appropriate map.
\item 
One can build examples of probability measures $\mu$ of positive dimension on the $d$-torus which do have escape of mass
on average (in the context of Theorems~\ref{theorem 1'},\ref{theorem 2'}). 
The constructions we suggest depend on the dimension.
 
For $d=1$ consider the set $D\subset [0,1)$ consisting of numbers with diverging c.f.e.\ coefficients.
Any probability measures supported in $D$ will produce an example of a measure with full escape of mass. As the dimension of $D$ is $\frac12$ (see~\cite{Good} or~\cite{Cheung1}), we conclude that such measures with positive dimension exist.

For $d=2$ one can take, in a similar manner, a probability measure of positive dimension supported on the set of singular vectors (recall that a vector $v\in\bR^d$ is said to be singular if the orbit $a_tu_v\PGL_{d+1}(\bZ)$ goes to infinity as $t\to\infty$). In ~\cite{Cheung1} the dimension of this set was calculated. For such measures there will be full escape of mass.

We do not know however if the further invariance assumption on the measure $\mu$ which appears in the statements of Theorems~\ref{theorem 1'},\ref{theorem 2'} actually excludes the possibility of escape of mass on average or even of escape of mass. 
\item We expect that the analogues for Theorems ~\ref{theorem 2},~\ref{theorem 2'} for higher dimensional torus still hold with some assumptions on the automorphism $\ga$ (or even if $\ga$ is an epimorphism). Using the high entropy method developed in ~\cite{EK}, it can be proved (and will be done elsewhere), that for any $d\ge3$, if $\ga$ is an automorphism of the $d$-torus with characteristic polynomial having only real roots which are distinct in absolute value, and $\mu$ is a $\ga$-invariant and ergodic measure of dimension greater than $1$, then $\mu$ almost any point is WA. Moreover if $\mu$ is friendly, then $\mu$ almost any point is not DI and has property C. 
\item Theorems ~\ref{theorem 1'},~\ref{theorem 2'} seem to have many applications to Diophantine approximations and the list of properties in their statements is not complete. 
For example, the third named author proved that for $s\in\bR$, if ~\eqref{density1} holds, then 
the 2-dimensional lattice $u_s\PGL_2(\bZ)$ (see \S\S~\ref{space of lattices}) satisfies the generalized Littlewood conjecture. For a proof of this statement and the discussion on the generalized Littlewood conjecture, the reader
is referred to ~\cite{ShapGLC}. 
\item M. Boshernitzan reported to us that a stronger version of Theorem ~\ref{theorem 3} holds for the special case of  quadratic irrationals.
\item B. de Mathan and and O. Teuli\'e have conjectured in ~\cite{De} that for any prime $p$ and for any irrational number $x\in[0,1]$, if we denote by $\tilde{c}(x)=\sup_n a_n(x)$, 
(where $a_n(x)$ are the coefficients in the continued fraction expansion of $x$), one has $\sup_\ell \tilde{c}(p^\ell x)=\infty$. In \cite{E-Kl} it was shown that the set of exceptions to de Mathan-Teuli\'e's conjecture
is of Hausdorff dimension zero. Although in Theorem ~\ref{theorem 3} we allow to multiply $x$ by a much bigger set
of integers than powers of a single prime, our result does not follow from de Mathan-Teuli\'e's conjecture because of the 
fundamental difference between the definitions of $c(x)$ using $\limsup$ and $\tilde{c}(x)$ using $\sup$. 
\end{enumerate}
\end{remarks}
\begin{acknowledgments}
We would like to express our gratitude to Barak Weiss and Dmitry Kleinbock for valuable suggestions and remarks. We also acknowledge the Max Planck institute and the program ``Dynamical Numbers'', held there on June 2009, for their kind hospitality. The first named author is grateful to Elon Lindenstrauss for various discussions concerning Theorem ~\ref{theorem 3}.
The second named author would like to thank Michael Boshernitzan for suggesting
and discussing many questions, one of which is solved in this paper. We would also like to thank the referee for helpful remarks.
\end{acknowledgments}
\section{Preliminaries}\label{preliminaries}
Most of the arguments appearing in our proofs are dynamical. In this section we present the dynamical systems in which
our discussion takes place and give the necessary preliminaries needed to understand the proofs of the results stated in \S\S\ref{statements of results}.   
\subsection{Homogeneous spaces} 
Let $G$ be a second countable locally compact topological group and $\Ga<G$ a discrete subgroup. The space $G/\Ga$ is called a homogeneous space as $G$ acts transitively on $G/\Ga$ by left translation. The topology we take on $G/\Ga$ is the quotient topology which then makes the natural projection $G\to G/\Ga$ a covering map. When $G/\Ga$ supports
a $G$-invariant probability measure we say that $\Ga$ is a lattice in $G$. In this case, this probability measure
is unique and is denoted by $\mu_G$. If $\Ga<G$ is a lattice, then the support of $\mu_G$ equals of course $G/\Ga$. This simple fact is used without reference in our arguments.  In this paper we will be interested in a very restrictive family of examples. We now describe the most important one. 
\subsection{The space of lattices}\label{space of lattices}
Fix $d\ge 1$ and let $X=\PGL_{d+1}(\bR)/\PGL_{d+1}(\bZ)$. It is well known that $\PGL_{d+1}(\bZ)<\PGL_{d+1}(\bR)$ is a lattice. The space $X$ can be identified with
the space of unimodular lattices in $\bR^{d+1}$ (i.e.\ of covolume 1) in the following manner: Given a coset $g\PGL_{d+1}(\bZ)$ we choose a matrix in $\GL_{d+1}(\bR)$ representing $g$ and denote it also by $g$. We then take the lattice spanned by the columns of $g$ and normalize it to have covolume 1. The reader should check that this defines a bijection between $X$ and the space of unimodular lattices in $\bR^{d+1}.$ The group $\SL_{d+1}(\bR)$ is
mapped in a natural way into $\PGL_{d+1}(\bR)$ and hence acts on $X$ by left translation. When we think of points of $X$ as lattices in $\bR^{d+1}$, this action translates to the linear action of $\SL_{d+1}(\bR)$ on $\bR^{d+1}.$ The following is known as Mahler's compactness criterion. It gives a geometric criterion for divergence in $X$ and in particular, shows that $X$ is not compact:
\begin{theorem}[Mahler's compactness criterion]\label{Mahler}
A subset $C\subset X$ is bounded (i.e.\ its closure is compact) if and only if there exists $\eps>0$ such that 
for any lattice $\Lam\in C$, $\Lam\cap B_\eps(0)=\set{0}$ i.e.\ if and only if there exists a uniform lower bound for
the lengths of nonzero vectors belongings to points in $C$.
\end{theorem}
We denote for $t\in \bR$ and $v\in\bR^d$, 
\begin{equation}\label{notation 1}
a_t=\operatorname{diag}(e^t,\dots,e^t,e^{-dt}),\quad 
u_v=
\pa{
\begin{array}{ll}
I_d& -v\\
0& 1
\end{array}
}\in \PGL_{d+1}(\bR),
\end{equation}
where $I_d$ is the $d\times d$ identity matrix. The mysterious minus sign in front of $v$ in ~\eqref{notation 1} is explained in the discussion in Appendix ~\ref{appendix}. 
Note that $\set{a_t}_{t\in\bR},\set{u_v}_{v\in\bR^d}$ are subgroups of $\PGL_{d+1}(\bR)$.
In the base of our arguments lies the identification of the $d$-torus $\bR^d/\bZ^d$ with the periodic orbit 
of the group $\set{u_v}_{v\in\bR^d}$ through the identity coset, 
\begin{equation}\label{d-torus as periodic orbit}
\textrm{for all } v\in\bR^d, \quad v+\bZ^d\leftrightarrow u_v\PGL_{d+1}(\bZ).
\end{equation}
Using this identification, many of the Diophantine properties of a vector $v\in\bR^d$, correspond to dynamical properties of the orbit $\set{a_tu_v\PGL_{d+1}(\bZ)}_{t>0}$. This is the content of Lemmas ~\ref{WA vectors and dynamics}
-- ~\ref{all patterns}. Although these  are probably well known, the proofs of Lemmas ~\ref{WA vectors and dynamics},~\ref{c.f.e} and ~\ref{all patterns} appear in Appendix ~\ref{appendix} for completeness of the exposition. The following lemma is essentially contained in Theorem 2.20 from \cite{Dani85}: 
\begin{lemma}\label{WA vectors and dynamics}
For any $\eps>0$ there exists a compact set $K_\eps\subset \PGL_{d+1}(\bR)/\PGL_{d+1}(\bZ)$ such that for any 
$v\in\bR^d$, if the inequality $\norm{v-\frac{\vec{m}}{n}}_\infty<\frac{\eps}{n^{1+1/d}}$ has only finitely many 
solutions $\vec{m}\in\bZ^d, n\in\bN$, then for large enough $T$, $a_tu_v\PGL_{d+1}(\bZ)\in K_\eps$ for $t>T$.
In particular if the vector $v\in\bR^d$ is not WA then the orbit $\set{a_tu_v\PGL_{d+1}(\bZ)}_{t\ge 0}$ is bounded.
\end{lemma}
\begin{lemma}\label{density imply}
Let $d\ge 2$. If $v\in\bR^d$ is such that 
$$\overline{\set{a_tu_v\PGL_{d+1}(\bZ)}}_{t>0}= \PGL_{d+1}(\bR)/\PGL_{d+1}(\bZ),$$  
then $v$ is not DI and has property C.
\end{lemma}
\begin{proof}
The proof of this lemma follows from Corollaries 4.6,4.8 in ~\cite{Shap1} and Proposition 2.1 in ~\cite{KW}.
\end{proof} 
The following lemma is left as an exercise
\begin{lemma}\label{WA invariant}
The class of WA points in the $d$-torus is invariant under the natural action of $M_d(\bZ)\cap \GL_d(\bQ)$.
\end{lemma}
\subsection{Connection with continued fraction expansion} We identify the circle $\bR/\bZ$ with the interval $[0,1)$. For each irrational $s\in [0,1)$, there exists a unique infinite sequence of positive integers $a_n(s)=a_n$ such that the sequence 
\begin{equation}\label{q.f}
[a_1,\dots,a_n]=\frac{1}{a_1+\frac{1}{a_2+\frac{1}{\ddots\frac{1}{a_n}}}}
\end{equation}
converges to $s$. This correspondence is a homeomorphism between $\bN^\bN$ and the irrational points on the circle.
We then denote $s=[a_1,a_2,\dots]$ and refer to the sequence $a_n(s)$ as the continued fraction expansion (abbreviated c.f.e.)\ of $s$. We denote, as in Theorem ~\ref{theorem 3} $c(s)=\limsup a_n(s)$. 
\begin{lemma}\label{c.f.e}
For any $N\in \bN$, there exists a compact set $K_N\subset  \PGL_2(\bR)/\PGL_2(\bZ)$ such that if $s\in \bR/\bZ$ is irrational and 
 $c(s)<N$ then the orbit $\set{a_tu_s\PSL_2(\bZ)}_{t\ge T}$ is contained in $K_N$ for large enough $T$ (which depends of course on $s$). 
\end{lemma}
We say that the c.f.e.\ of an irrational $s\in\bR/\bZ$ contains all patterns if given a finite sequence of integers 
$b_1,\dots,b_n$, there exists $k$ such that $a_{k+i}(s)=b_i$ for any $1\le i\le n$.
\begin{lemma}\label{all patterns}
If $s\in\bR/\bZ$ is such that $\overline{\set{a_tu_s\PGL_2(\bZ)}}_{t>0}=\PGL_2(\bR)/\PGL_2(\bZ)$, then the c.f.e.\ of $s$ 
contains all patterns.
\end{lemma}
%
%
%
%
%
%
%
%
%
%
\subsection{Escape of mass} 
Given a probability measure $\mu$ on $\bR^d/\bZ^d$, we may think of it (see ~\eqref{d-torus as periodic orbit}) as a measure supported on the periodic orbit $$\set{u_v\PGL_{d+1}(\bZ)}_{v\in\bR^d}\subset \PGL_{d+1}(\bR)/\PGL_{d+1}(\bZ).$$ This enables us to define
\begin{definition}\label{no escape in avg}
We say  that $\mu$ has no escape of mass on average with respect to $\set{a_t}_{t\ge 0}$ if any weak$^*$ limit of 
$\frac{1}{T}\int_0^T (a_t)_*\mu dt$ is a probability measure on $\PGL_{d+1}(\bR)/\PGL_{d+1}(\bZ).$
\end{definition}
We can now state Theorem 5.3 from ~\cite{R} which will be needed to prove Theorems ~\ref{theorem 1},~\ref{theorem 1'}.
\begin{theorem}[Theorem 5.3 from ~\cite{R}]\label{Ronggang's thm}
Let $\mu$ be a probability measure on $\bR^d/\bZ^d$ of dimension $\ka$ such that $\mu$ has no escape of mass on average with respect to $\set{a_t}_{t\ge0}$. Then any weak$^*$ limit, $\nu$, of $\frac{1}{T}\int_0^T(a_t)_*\mu dt$
satisfies $h_\nu(a_1)\ge (d+1)\ka$. In particular, if $\ka>0$ then $h_\nu(a_1)>0.$
\end{theorem}
\subsection{Group action on measures}\label{acting on measures}
Let $X= G/\Ga$ be a homogeneous space ($G$ a locally compact group and $\Ga$ a discrete subgroup of $G$). $G$ acts on $X$ by left translations. This action induces an action of $G$ on the space of Borel probability measures on $X$. Given a probability measure $\mu$ on $X$ and $g\in G$, we denote by $g_*\mu$ the probability measure defined by the equation
\begin{equation}\label{left action on prob}
 \int_X f(x)dg_*\mu(x)=\int_X f(gx)d\mu(x)
\end{equation} 
for any  $f\in C_c(X)$. $\mu$ is said to be $g$-invariant if $g_*\mu=\mu$. 
Given a subgroup $H<G$,
the set of $H$-invariant probability measures will be denoted by $\cM_X(H)$.\\
Let $H<G$ be a commutative closed group and let $\mu\in\cM_X(H)$. The ergodic decomposition of $\mu$ with respect to $H$ is the unique Borel probability measure $\theta_H$ concentrated on the extreme points of $\cM_X(H)$ (i.e.\ the extreme points have $\theta_H$-measure 1) and having $\mu$ as its center of mass. Existence and uniqueness of the ergodic decomposition follow from Choquet's theorem. We say that an ergodic $H$-invariant measure $\mu_0$, appears as a component with positive weight in the ergodic decomposition of $\mu$ with respect to $H$, if $\theta_H(\set{\mu_0})>0$.
An equivalent (and perhaps simpler) condition is the existence of a constant $c>0$, such that for any nonnegative function $f\in C_c(X)$ one has
$\int_Xfd\mu\ge c\int_X fd\mu_0$.  \\
Let $H'<H$ be a closed subgroup. 
If $\mu_0$ is ergodic with respect to $H'$ (and hence with respect to $H$), then it appears with positive weight in the ergodic decomposition of $\mu$ with respect to $H$, if and only if it appears as a component with positive weight in the ergodic decomposition with respect to $H'$.\\
$H$ acts on $\cM_X(H')$ and as $H'$ acts trivially, this action induces an action of the quotient $H/H'$ on $\cM_X(H')$. Denote the natural projection from $H$ to $H/H'$ by $g\mapsto\hat{g}$. Let $\mu\in\cM_X(H')$. If the quotient $H/H'$ is compact, one can define an $H$-invariant probability measure 
$$\tilde{\mu}=\int_{H/H'}\hat{g}_*\mu d\hat{g},$$
where $d\hat{g}$ is the Haar probability  measure on $H/H'$. The meaning of this equation is that 
$$\int_X f(x)d\tilde{\mu}=\int_{H/H'}\pa{\int_X f(x)d\hat{g}_*\mu}d\hat{g}$$
for any  $f\in C_c(X)$.
For $b\in H'$ and $g\in H$, the entropies $h_\mu(b),h_{g_*\mu}(b)$ are equal. This implies that
$h_{\tilde{\mu}(b)}=h_\mu(b)$ too.
We shall need the following theorem about entropy (see ~\cite{EL1} for the proof). 
\begin{theorem}[Upper semi continuity of entropy]\label{u.s.c of entropy}
Let $X=G/\Ga$ be as above and let $b\in G$. Let $\mu_n$ be a sequence of 
$b$-invariant probability measures converging in the weak$^*$ topology to a probability measure $\mu$ (which is automatically $b$-invariant). Then $h_\mu(b)\ge\limsup h_{\mu_n}(b).$
\end{theorem} 
\section{Proofs of Theorems ~\ref{theorem 2},~\ref{theorem 2'}}
In this section $G=\PGL_3(\bR),\Ga=\PGL_3(\bZ)$ and $X= G/\Ga$. The identity coset in $X$ will be denoted by 
$\bar{e}$. We use the notation of ~\eqref{notation 1} and the identification of ~\eqref{d-torus as periodic orbit}.
Hence the 2-torus $ \bR^2/\bZ^2$ is identified with the periodic orbit $\set{u_v\bar{e}:v\in\bR^2}$ of the two dimensional unipotent group $\set{u_v}_{v\in\bR^2}<G$. This enables us to view the measure $\mu$ from the statement of Theorems ~\ref{theorem 2},~\ref{theorem 2'}, as a measure supported on this periodic orbit. Let $\ga$ be as in the statement of Theorem ~\ref{theorem 2}. Under this identification,
the action of $\ga$ translates to the action from the left of 
\begin{equation}\label{action on the torus}
\ga'=\pa{
\begin{array}{ll}
\ga&0\\
0& 1
\end{array}
}\in \Ga.
\end{equation}
Hence, the assumptions of Theorem ~\ref{theorem 2} translate to $\mu$ being of positive dimension, $\ga'$-invariant, and ergodic.
\begin{proof}[Proof of Theorem ~\ref{theorem 2}]
Assume that the statement of the theorem is false. 
As the set of WA points on the torus is $\ga'$-invariant (see Lemma ~\ref{WA invariant}), 
it follows from ergodicity and Lemma ~\ref{WA vectors and dynamics}, that for $\mu$-almost any $x\in X$, the orbit $\set{a_tx}_{t\ge 0}$, is bounded. 
Let $K_i$ be an increasing sequence of compact subsets exhausting $X$.
We shall build an invariant measure on $X$ having the Haar measure appearing as a component with positive weight in its ergodic decomposition, while at the same time this measure will
be the sum of invariant measures supported on the sets $K_i$. This contradicts the uniqueness of the ergodic decomposition.

To this end we define 
$$E_i=\set{x\in \operatorname{supp}(\mu) : \set{a_tx}_{t\ge 0}\textrm{ is contained in $K_i$ but not in $K_{i-1}$}}.$$
Hence, $E_i$ form a partition (up to a null set) of the support of $\mu$. Denote by $\mu_i$ the restriction of $\mu$ to $E_i$. Hence $\mu=\sum\mu_i$. We denote $\mu_i^T=\frac{1}{T}\int_0^T (a_t)_*\mu_idt$ and $\mu^T=\sum\mu_i^T$. Let $T_j\to\infty$ be chosen so that the sequences $\mu_i^{T_j},\mu^{T_j}$ converge weak$^*$ to some measures $\nu_i,\nu$ respectively. Since for any $t\ge 0$, $a_t(E_i)\subset K_i$, $\nu_i$ is supported in $K_i$ and there could be no escape of mass and $\nu$ is a probability measure.
$\nu$ and the $\nu_i$'s are $a_t$-invariant and $\nu=\sum\nu_i$. In particular, the ergodic decomposition of $\nu$ with respect to $\set{a_t}_{t\in\bR}$ is the sum of the ergodic decompositions of the $\nu_i$'s. As $\nu_i$ is supported in $K_i$, we deduce that the $G$-invariant probability measure $\mu_G$, cannot appear as a component with positive weight, in the ergodic decomposition of $\nu$ with respect to $\set{a_t}_{t\in\bR}.$ 
Since the action of $\ga'$ commutes with $a_t$, $\mu^{T_j}$ is $\ga'$-invariant for any $j$ and as a consequence $\nu$ is $\ga'$-invariant too. Note also that for any $T$ we have the following equality of entropies: $h_\mu(\ga')=h_{\mu^T}(\ga')$. Hence it follows from Remark~\ref{dimension is entropy} and Theorem ~\ref{u.s.c of entropy} that $h_\nu(\ga')>0$.
From our assumption on the hyperbolicity of $\ga$ (which in this case implies $\bR$-diagonability), it follows that the group, $H$, generated by $\set{a_t}_{t\in\bR}$ and $\ga'$, is cocompact in a maximal $\bR$-split torus $T$ in $G$. 
The desired contradiction now follows from Corollary ~\ref{positive entropy 4} below, which in turn follows 
from  the following theorem from ~\cite{EKL}. 
\end{proof}
\begin{theorem}[Theorem 1.3 from ~\cite{EKL}]\label{positive entropy 3}
Let $\nu$ be a Borel probability measure on $X= \PGL_3(\bR)/\PGL_3(\bZ)$ which is invariant under the action of a maximal $\bR$-split torus $T<G=\PGL_3(\bR).$ If there exists $b\in T$ which acts with positive entropy with respect to $\nu$, then the $G$-invariant probability measure $\mu_G$, appears as a component with positive weight in the ergodic decomposition of $\nu$ with respect to $T$.
\end{theorem} 
\begin{corollary}\label{positive entropy 4}
Let $\nu$ be a Borel probability measure on $X$ which is invariant under the action of a group $H$ which is cocompact in a maximal $\bR$-split torus $T<G.$ If there exists $b\in H$ which acts with positive entropy with respect to $\nu$, then the $G$-invariant probability measure $\mu_G$, appears as a component with positive weight in the ergodic decomposition of $\nu$ with respect to $H$.
\end{corollary}
\begin{proof}  
Denote the natural projection $T\to T/H$ by $g\mapsto \hat{g}$. Define 
$$\lam=\int_{T/H}\hat{g}_*\nu d\hat{g},$$
where $d\hat{g}$ is the Haar measure in $T/H$ (recall the discussion of \S\S~\ref{acting on measures}). 
$\lam$ is a $T$-invariant measure on $X$ and $h_\lam(b)=h_\nu(b)$, hence Theorem ~\ref{positive entropy 3} implies that $\mu_{G}$ appears with positive weight in the ergodic decomposition of $\lam$ with respect to $T$. By the Howe-Moore theorem $\mu_G$ is $H$-ergodic, hence we conclude that $\mu_{G}$ appears with positive weight in the ergodic decomposition of $\lam$ with respect to $H$. The ergodic decomposition of $\nu$ with respect to $H$ is a probability measure $\theta$, supported on the extreme points of $\cM_X(H)$, which is the set of $H$-invariant probability measures on $Y$,
having $\nu$ as its center of mass. The ergodic decomposition of $\lam$ with respect to $H$ is $\theta'=\int_{T/H}\hat{g}\theta d\hat{g}.$ This equation means that $\theta'$ is the probability measure on $\cM_X(H)$, characterized by the following equation:
\begin{equation}\label{positive weight 1}
\int_{\cM_X(H)} F(\vphi) d\theta'(\vphi)=\int_{T/H}\int_{\cM_X(H)}F(\hat{g}_*\vphi)d\theta(\vphi)d\hat{g}
\end{equation}
for any  $F\in C(\cM_X(H))$. 
In order to show that $\theta'(\set{\mu_{G}})>0$ and conclude the proof, we need to show that for any open neighborhood $\mu_{G}\in V\subset \cM_X(H)$, $\theta'(V)>\al$ for some positive constant $\al$. Let $V$ be such an open neighborhood. Let $U\subset V$ be another open neighborhood of $\mu_{G}$ such that there exists a bump function $F$ which equals 1 on $U$ and vanishes outside $V$. Let $U'\subset U$ be a smaller neighborhood of $\mu_G$, such that
\begin{equation}\label{positive weight 3}
U'\subset \cap_{\hat{g}\in T/H} \hat{g}_*(U).
\end{equation}
The existence of $U'$ follows from the compactness of $T/H$ and the $G$-invariance of $\mu_G$. 
Then 
\begin{equation}\label{positive weight 2}
\theta'(V)\ge \int_{\cM_X(H)} Fd\theta'=\int_{T/H}\int_{\cM_X(H)}F(\hat{g}_*\vphi)d\theta(\vphi)\hat{g}.
\end{equation}
By construction, for any $\hat{g}\in T/H$
\begin{equation}\label{positive wight 4}
\int_{\cM_X(H)}F(\hat{g}_*\vphi)d\theta(\vphi)\ge\theta(U')\ge \al,
\end{equation}
where $\al=\theta(\set{\mu_{G}})$ is positive by assumption.
\end{proof}
\begin{proof}[Proof of Theorem ~\ref{theorem 2'}]
As $\ga'$ and $a_t$ commute, the set $F=\set{x\in X : \overline{\set{a_tx}}_{t>0}\ne X}$ is $\ga'$-invariant. Assume to get a contradiction that $\mu(F)>0$. It follows from the ergodicity that $\mu(F)=1$. Let $\set{U_i}$ be
a countable base for the topology of $X$. Define recursively 
\begin{align*}
E_1&=\set{x\in \operatorname{supp}(\mu) : \set{a_tx}_{t>0}\cap U_1=\emptyset}\mbox{ and}\\
E_n&=\set{x\in \operatorname{supp}(\mu) : \set{a_tx}_{t>0}\cap U_n=\emptyset}\setminus E_{n-1}
\end{align*}
for any $n>1$.
Hence, $\set{E_i}$ form a partition up to a null set of the support of $\mu$. We continue as in the proof of Theorem ~\ref{theorem 2} using the same notation as there. We now highlight the differences between the arguments: In the proof of Theorem ~\ref{theorem 2} we used the fact that $\mu_i^T$ is compactly supported in order to pass to a weak$^*$ limit without losing mass. Here we do not know that $\mu_i^T$ is compactly 
supported, instead we use our further assumption that any weak$^*$ limit of $\mu^T$ is a probability measure. Another difference is that in the proof of Theorem ~\ref{theorem 2}, $\nu_i$ was supported in a compact set and hence could not have $\mu_G$ appear as a component with positive weight in its ergodic decomposition with respect to 
$\set{a_t}_{t\in\bR}$. Here the reason that $\nu_i$ cannot have $\mu_G$ appearing as a component with
positive weight is that $\nu_i(U_i)=0.$ 

To end the proof we note that Lemmas ~\ref{WA vectors and dynamics}, ~\ref{density imply} imply that for $\mu$-almost any $v\in\bR^2/\bZ^2$, $v$ is WA, not DI, and has property C.
\end{proof} 
\section{Proof of Theorems ~\ref{theorem 1},~\ref{theorem 1'}}
\subsection{Preparations}
 Let $\bG=\PGL_2$ and $S=\set{p_1,\dots,p_k,\infty}$, where the $p_i$'s are the primes appearing in the prime decomposition of the number $n$ appearing 
in the statement of Theorem ~\ref{theorem 1}. We denote 
 \begin{equation}\label{groups}
G_\infty=\bG(\bR),\quad G_f=\prod_1^k\bG(\bQ_{p_i}),\quad G_S=G_\infty\times G_f, \quad K=\prod_1^k\bG(\bZ_{p_i}). 
\end{equation}
Denote $\Ga_S=\bG(\bZ[\frac{1}{p_1}\dots\frac{1}{p_k}])$ and $\Ga_\infty=\bG(\bZ)$.
 We shall
abuse notation (as usual) and identify $\Ga_S$ with its various diagonal embeddings in $G_S,G_f$ etc. The meaning should be clear from the context.  
$\Ga_S,\Ga_\infty$ are lattices in $G_S,G_\infty$ respectively. Nonetheless, $\Ga_S$ is dense in $G_f$.  Let $X= G_\infty/\Ga_\infty$ and
 $Y= G_S/\Ga_S$. 
 We denote the identity cosets in both spaces by $\bar{e}$. 
 The elements of $G_S$ will be denoted by $(g_\infty,g_f)$ where $g_\infty\in G_\infty$ and $g_f\in G_f$.
Denote by
\begin{equation}\label{pi}
\pi: Y\to K\backslash Y=K\backslash G_S/\Ga_S,
\end{equation}
the natural projection. The double coset space $K\backslash G_S/\Ga_S$ can be identified with $X$ in the following manner: 
Given a double coset $K(g_\infty,g_f)\Ga_S$ 
as $K$ is an open subgroup of $G_f$ and $\Ga_S< G_f$ is dense, there exists $\ga\in\Ga_S$ such that $g_f\ga\in K$. We 
then identify $K(g_\infty,g_f)\Ga_S$ with $g_\infty\ga\bar{e}\in X$. The reader should check that this map is indeed
well defined, a bijection, and respects the topologies. In other words the map $\pi:Y\to X$ is defined by
$$\pi\pa{(g_\infty,g_f)\bar{e}}= g_\infty\bar{e}\textrm{ if }g_f\in K.$$
 $G_S,G_\infty$  act on $Y,X$ by left translation respectively. The action of $G_\infty$ on $X$ is via $\pi$ a factor of the action of $G_\infty\times\set{e_f}$ on $Y$.  
\subsection{Proofs}
\begin{proof}[Proof of Theorem ~\ref{theorem 1}]
We identify $ \bR/\bZ$ with the  periodic orbit of the horocycle flow $\set{u_t}_{t\in\bR}$ through $\bar{e}\in X$ (see ~\eqref{notation 1},~\eqref{d-torus as periodic orbit}).
Under this identification, the map $\times n$ becomes the map $u_s\bar{e}\mapsto u_{ns}\bar{e}$. This identification enables us to view the measure $\mu$ from Theorem ~\ref{theorem 1} 
as a probability measure supported on this periodic orbit. The next thing we wish to do is to lift this measure to a measure on $Y$. We do so by pushing it with the map $u_t\bar{e}\mapsto (u_t,e_f)\bar{e}$ defined for $t\in[0,1)$. We denote the resulting measure on $Y$ by $\nu_1$. It is obvious that $\pi_*(\nu_1)=\mu$. 
We let $b=\operatorname{diag}(n,1)\in \Ga_S$ and note that
the action of $b$ on $Y$, when restricted to $\set{(u_s,e_f)\bar{e}:s\in\bR}$, factors via $\pi$ to the
map $\times n$ on the circle; i.e.\ the following diagram commutes:
\begin{equation*}
\xymatrix{
(u_s,e_f)\bar{e}\ar[r]^b\ar[d]_\pi &(u_{ns},e_f)\bar{e}\ar[d]^\pi\\
u_s\bar{e}\ar[r]^{\times n} & u_{ns}\bar{e}
}
\end{equation*}
Although $\mu$ is $\times n$-invariant, $\nu_1$ is not invariant under the action of 
$b$ on $Y$. We replace it by a different measure which is invariant under $b$ and projects to $\mu$ by the following 
procedure: We denote $\nu_N=\frac{1}{N}\sum_{i=0}^{N-1}b^i_*(\nu_1)$. Note that for any $N$, $\pi_*(\nu_N)=\mu$. Let 
$\nu$ be a weak$^*$ limit of the sequence $\nu_N$. It follows that $\pi_*(\nu)=\mu$ and in particular, 
that $\nu$ is a probability measure (note that here we used the fact that the fibers of $\pi$ are compact)\footnote{One could modify the above construction and first lift the measure $\mu$ from $\bR/\bZ$ to $\bR\times\bQ_S/\bZ_S$ and then average the lift to get invariance under the (invertible extension of) $\times n$ and only then identify the resulting measure $\nu$ with a measure on $Y$ which projects to $\mu$.}.
To summarize  what we established so far, we constructed a $b$-invariant probability measure, $\nu$, on $Y$ such that  $\pi:(Y,\nu,b)\to(X,\mu,\times n)$, is a factor map. 
Assume that the statement of Theorem ~\ref{theorem 1} is false. It follows from Lemmas ~\ref{WA invariant} and  ~\ref{WA vectors and dynamics} that for $\mu$-almost any $x\in X$, $\set{a_tx:t\ge 0}$ is bounded. Let $K_i$ be an increasing sequence of compact subsets exhausting $X$. Let $$E_i=\set{x\in \operatorname{supp}(\mu) : \set{a_tx}_{t\ge 0} \textrm{ is contained in } K_i \textrm{, but not in } K_{i-1}}.$$ 
Thus $E_i$ form a partition (up to a null set) of the support of $\mu$. Denote by $\mu_i$ the restriction of $\mu$ to $E_i$, hence $\mu=\sum_i\mu_i$.
We denote for $T>0$, $\nu^T=\frac{1}{T}\int_0^T(a_t,e_f)_*(\nu)dt.$  Then 
\begin{equation}\label{factors}
\pi_*(\nu^T)=\frac{1}{T}\int_0^T (a_t)_*(\mu)dt=\sum_i\frac{1}{T}\int_0^T (a_t)_*(\mu_i)dt.
\end{equation}
 Denote $\mu_i^T=\frac{1}{T}\int_0^T (a_t)_*(\mu_i)dt$ and $\mu^T=\sum_i\mu_i^T$. Thus, ~\eqref{factors} becomes 
 $\pi_*(\nu^T)=\mu^T=\sum\mu_i^T$. Let $T_j\to\infty$ be chosen such that all the following sequences converge in the
 weak$^*$ topology: $\nu^{T_j},\mu_i^{T_j},\mu^{T_j}$. Denote their corresponding limits by 
 $\tilde{\nu},\tilde{\mu}_i,\tilde{\mu}$ respectively. It is evident that $\tilde{\nu}$ is $(a_t,e_f)$-invariant, 
 while $\tilde{\mu},\tilde{\mu}_i$ are $a_t$-invariant. As the fibers of $\pi$ are compact, we can deduce that 
 $\pi_*(\tilde{\nu})=\tilde{\mu}=\sum_i\tilde{\mu}_i$. Moreover since 
$\tilde{\mu}_i$ is supported in $K_i$, there is no escape of mass and $\tilde{\nu},\tilde{\mu}$ are probability measures. 
We will derive the desired contradiction by using the following lemma:
\begin{lemma}\label{Haar component}
In the ergodic decomposition of $\tilde{\mu}$ with respect to $\set{a_t}_{t\in\bR}$, the $G_\infty$-invariant measure $\mu_{G_\infty}$ has positive weight.
\end{lemma}
To finish the proof of the theorem, note that since for each $i$, $\tilde{\mu}_i$ is $a_t$-invariant, the ergodic decomposition of $\tilde{\mu}$ with respect to the action of $\set{a_t}_{t\in\bR}$ is the sum of the corresponding ergodic decompositions of the $\tilde{\mu}_i$'s which are supported in $K_i$ and hence cannot have $\mu_{G_\infty}$ appearing with positive weight in their ergodic decomposition.
\end{proof}
In the proof of Lemma ~\ref{Haar component} we will use the following simplification of Theorem 1.1 from ~\cite{L}. To state it we use the notation
from the beginning of this subsection and we denote by $\bT$ the subgroup of $\bG$ consisting of diagonal matrices.
\begin{theorem}\label{positive entropy 1}
Let $\tilde{\nu}$ be a probability measure on $Y$ which is invariant under the action of $\bT(\bR)$, has positive entropy with respect some (hence any) element in $\bT(\bR)$ and is invariant under the action of of a noncompact subgroup of $G_f$. 
Then in the ergodic decomposition of $\pi_*(\tilde{\nu})$ with respect to $\bT(\bR)$, the $G_\infty$-invariant measure $\mu_{G_\infty}$ appears as a component with positive weight.
\end{theorem}  
\begin{proof}[Proof of Lemma ~\ref{Haar component}]
By construction, the measure $\tilde{\nu}$ is invariant under the group generated by $\bT(\bR)=\set{(a_t,e)}_{t\in\bR}$ and $(b,b)$ 
(here we use the fact that $(a_t,e)$ and $(b,b)$ commute). In particular $\tilde{\nu}$ is invariant under a noncompact subgroup of $G_f$.
It follows from the positivity of the dimension of $\mu$ and Theorem~\ref{Ronggang's thm}, that 
$h_{\tilde{\mu}}(a_1)>0$. Then, since $(X,\tilde{\mu},a_t)$ is a factor of $(Y,\tilde{\nu},(a_t,e_f))$, we must have
$h_{\tilde{\nu}}((a_1,e_f))>0$. We see that the conditions of Theorem ~\ref{positive entropy 1} are satisfied and 
as a consequence that $\tilde{\nu}=\pi_*(\tilde{\nu})$ has $\mu_{G_\infty}$ appearing as a component with positive weight in the ergodic decomposition of it with respect to the action of $\set{a_t}_{t\in\bR}$ as desired.
\end{proof}
In order to complete the proof of Theorem ~\ref{theorem 1'} we shall need the following lemma:
\begin{lemma}\label{invariance}
The set $F=\set{s\in\bR/\bZ:\overline{\set{a_tu_s\Ga_\infty}}_{t\ge 0}\ne X}$ is $\times n$-invariant.
\end{lemma}
\begin{proof}[Proof of Theorem ~\ref{theorem 1'}]
Let $F$ be as in Lemma ~\ref{invariance}. Assume to get a contradiction that $\mu(F)>0$. It follows from the ergodicity that $\mu(F)=1$. Let $\set{U_i}$ be
a countable base for the topology of $X$. Define recursively 
\begin{align*}
E_1&=\set{x\in \operatorname{supp}(\mu) : \set{a_tx}_{t>0}\cap U_1=\emptyset}\mbox{ and}\\
E_n&=\set{x\in \operatorname{supp}(\mu) : \set{a_tx}_{t>0}\cap U_n=\emptyset}\setminus E_{n-1}.
\end{align*}
for any $n>1$.
Hence, $\set{E_i}$ form a partition up to a null set of the support of $\mu$. Denote by $\mu_i$, the restriction
of $\mu$ to $E_i$. We continue as in the proof of Theorem ~\ref{theorem 1} using the same notation as there. We now highlight the differences between the arguments: In the proof of Theorem ~\ref{theorem 1} we used the fact that $\mu_i^T$ is compactly supported in order to pass to a weak$^*$ limit without losing mass. Here we do not know that $\mu_i^T$ is compactly 
supported, instead we use our further assumption that any weak$^*$ limit of $\mu^T$ is a probability measure. In particular $\tilde{\mu}(X)=1$. This in turn implies that $\tilde{\nu}(Y)=1$. Another difference is that in the proof of Theorem ~\ref{theorem 1}, $\mu_i$ was supported in a compact set and hence could not have $\mu_{G_\infty}$ appear as a component with positive weight in its ergodic decomposition with respect to 
$\set{a_t}_{t\in\bR}$. Here the reason that $\mu_i$ cannot have $\mu_{G_\infty}$ appearing as a component with
positive weight is that $\mu_i(U_i)=0.$

After establishing the density of $\set{a_tu_s\Ga_\infty}_{t\ge 0}$ for $\mu$-almost any $s\in\bR/\bZ$, Lemma ~\ref{all patterns} implies that the continued fraction expansion of any such $s$ contains any given pattern.
\end{proof}
In order to prove Lemma ~\ref{invariance} we shall need the following lemma which follows immediately from ergodicity of the $a_t$ action on $X$:
\begin{lemma}\label{nowhere dense}
Let $C\subset X$ be closed and $\set{a_t}_{t\in\bR}$-invariant. Then either $C=X$ or $C$ has empty interior.
\end{lemma}
\begin{proof}[Proof of Lemma ~\ref{invariance}]
Let us change notation and set 
$$G=\PGL_2(\bR),\quad\Ga_1=\PGL_2(\bZ),\quad\Ga_2=\operatorname{diag}(n^{-1},1)\Ga_1\operatorname{diag}(n,1),
\textrm{ and }\Ga=\Ga_1\cap \Ga_2.$$
Note that $\Ga$ is of finite index in both of the $\Ga_i$'s. It means that the natural projections 
$p_i:G/\Ga\to G/\Ga_i$ are finite covers. As such, they satisfy:
\begin{equation}\label{cover}
\textrm{For any } M\subset G/\Ga, \quad p_i(\overline{M})=\overline{p_i(M)}. 
\end{equation}
Let now $s\in\bR/\bZ$ be such that $ns\notin F$ i.e.\ such that $\overline{\set{a_tu_{ns}\Ga_1}}_{t>0}=G/\Ga_1$. We need to show that $s\notin F$ i.e.\ that the same holds for $s$ instead of $ns$. Assume first that
\begin{equation}\label{Ga1}
\overline{\set{a_tu_s\Ga_2}}_{t>0}=G/\Ga_2.
\end{equation}
It follows from ~\eqref{cover} that 
$p_2\pa{\overline{\set{a_tu_s\Ga}}_{t>0}}=G/\Ga_2$ so $\overline{\set{a_tu_s\Ga}}_{t>0}$ must have non empty interior in $G/\Ga$ (by Baire's category theorem for example) and in turn $p_1\pa{\overline{\set{a_tu_s\Ga}}_{t>0}}=\overline{\set{a_tu_s\Ga_1}}_{t>0}$ has nonempty interior in $G/\Ga_1$. Lemma ~\ref{nowhere dense} now implies that $\overline{\set{a_tu_s\Ga_1}}_{t>0}=G/\Ga_1$ as desired. We now argue
the validity of ~\eqref{Ga1}.  The fact that $\overline{\set{a_tu_{ns}\Ga_1}}_{t>0}=G/\Ga_1$ is equivalent to the set 
$$\set{a_t\operatorname{diag}(n,1)u_s \operatorname{diag}(n^{-1},1)\ga:t>0,\ga\in\Ga_1}$$
being dense in $G$. As $a_t$ and $\operatorname{diag}(n,1)$ commute, this is the same as to say that the set 
$$\set{a_tu_s\operatorname{diag}(n^{-1},1)\ga \operatorname{diag}(n,1):t>0,\ga\in\Ga_1}=\set{a_tu_s\ga:t>0,\ga\in\Ga_2}$$
is dense in $G$, which is exactly ~\eqref{Ga1}.
\end{proof}
\section{Proof of Theorem ~\ref{theorem 3}}
In this section we use the following notation. Let $\bA$ denote the ring of adeles, $\bG=\PGL_2$, $G=\bG(\bR)$,
$G'=\bG(\bA)$, $\Ga=\bG(\bZ)$ and $\Ga'=\bG(\bQ)$. $\Ga'$ is a lattice in $G'$ when embedded diagonally.
We denote elements of $G'$ as $(g_\infty,g_2,g_3,g_5\dots)$ and will abbreviate and denote them simply as $(g_\infty,g_{f})$, where $g_f=(g_2,g_3\dots)$. Let $\bT<\bG$ be the the subgroup consisting of 
(classes of) diagonal matrices and denote
$T'=\bT(\bA),T=\bT(\bR)$. We denote $X=G/\Ga$ and $Y=G'/\Ga'$.
$\bar{e}$ will denote the identity coset in both spaces.
 Define $\pi:Y\to X$ in the following way:
 For a point $y\in Y$, we choose a representative $(g_\infty,g_{f})\in G'$, for which $g_{f}\in\bG(\prod_p\bZ_p)$, and define $\pi(y)=g_\infty\bar{e}$. 
$\pi$ is well defined, continuous and has compact fibers. We use the notation and identification of ~\eqref{notation 1},~\eqref{d-torus as periodic orbit} and identify $\bR/\bZ$ with the periodic orbit of the horocycle $u_t$, through the identity coset $\bar{e}\in X$. 
We shall also need the following theorem of E.\ Lindenstrauss and its corollary.
\begin{theorem}[Theorem 1.5 of~\cite{L-aue}]\label{unique ergodicity}
The action of the group, $T'$, of adelic points of the torus $\bT=\set{\operatorname{diag}(*,*)}<\PGL_2$ on $Y=\PGL_2(\bA)/\PGL_2(\bQ)$ is uniquely ergodic.
\end{theorem}

\begin{corollary}\label{unique ergodicity 2}
Let $H<T'$ be a cocompact subgroup. Then there are no compactly supported $H$-invariant measures on $Y$.
\end{corollary}
\begin{proof}
Assume by way of contradiction that $\nu$ is a compactly supported $H$-invariant measure on $Y$. Define 
$$\tilde{\nu}=\int_{T'/H}\hat{g}_*\nu d\hat{g}.$$
Then $\tilde{\nu}$ is $T'$-invariant and compactly supported (because $H$ is cocompact in $T'$). This contradicts Theorem ~\ref{unique ergodicity}.
\end{proof} 

\begin{proof}[Proof of Theorem ~\ref{theorem 3}]
For any prime $p$, let $b_p= \operatorname{diag}(p,1)\in\Ga'.$ We denote the diagonal embedding of $b_p$ in $G'$ by the same 
letter. Note that for any $s\in\bR$,  $b_p(u_s,e_{f})\bar{e}=(u_{ps},e_{f})\bar{e}$ and in particular,
if $n=p_1\dots p_k$, then $\pi\pa{b_{p_1}\dots b_{p_k}(u_s,e_{f})\bar{e}}=u_{ns}\bar{e}$.
Assume that the statement of the theorem is false.
Thus, by Lemma ~\ref{c.f.e} there exists a compact set $K\subset X$ and an irrational $s\in[0,1)$ such that for any $n=p_1\dots p_k$, for large enough $t$,
$$K\ni a_tu_{ns}\bar{e}=\pi\pa{(a_t,e_{f})b_{p_1}\dots b_{p_k}(u_s,e_{f})\bar{e}}.$$
Hence, if we denote $K'=\pi^{-1}(K)\subset Y$ then for fixed $b_{p_1},\dots,b_{p_k}$ and all sufficiently large ~$t$
\begin{equation}\label{supported in K'}
 (a_t,e_{f})b_{p_1}\dots b_{p_k}(u_s,e_{f})\bar{e}\in K'.
\end{equation}
Let $C<T'$ be the semigroup generated by 
$(a_1,e_{f})$ and the $b_p$'s and let $H$ be the group generated by $C$. $H$ is cocompact in $T'$.
To see this note that the compact set 
$$\set{\textbf{a}=\pa{\diag{e^t,e^{-t}},a_2,a_3\dots}:t\in[0,1],a_p\in\bT(\bZ_p)},$$
contains a
fundamental domain for $H$ in $T'$; this follows from the fact that for any element
 $\textbf{a}=(a_\infty,a_2,a_3,\dots)\in T'$, for almost all primes $p$, the matrix $a_p$ has entries in $\bZ_p$ (see Remark~\ref{non cocompact subgroups}).
  
Let $F_n$ be a F\o lner sequence for $C$ and define 
\begin{equation}\label{mu n}
\mu_n=\frac{1}{\av{F_n}}\sum_{g\in F_n}g_*\del_{(u_s,e_{f})\bar{e}},
\end{equation}
where $\del_{(u_s,e_{f})\bar{e}}$ is the Dirac measure centered  at the point $(u_s,e_{f})\bar{e}$. Let $\mu$ be a weak$^*$ limit of $\mu_n$. It is $H$-invariant. On the other hand, we claim that if the F\o lner sequence is chosen appropriately then by ~\eqref{supported in K'}, it is a probability measure supported in 
$K'$. This contradicts Corollary ~\ref{unique ergodicity 2}. We define $F_n$ inductively in the following manner: We first choose a F\o lner sequence, $F_n'$, for the semigroup $C'$ generated only by the $b_p$'s. Then for a fixed $n$, there is some $T_n$ such that for any $g\in F_n'$ and for any $t>T_n$, $(a_t,e_{f})g(u_s,e_{f})\bar{e}\in K'$. It follows that there exists an integer $m_n>T_n$, such that if we define 
\begin{equation}\label{F_n}
F_n=F_n'\cup\set{(a_1,e_{f})^k}_1^{m_n},
\end{equation}
then the weight that $\mu_n$ from ~\eqref{mu n} gives to $K'$ is greater than $1-1/n$.
\end{proof}
\begin{remark}\label{non cocompact subgroups}
It is tempting to replace in the above argument the group $H$ by  the group generated by $(a_1,e_f)$ and 
the elements $b_p^k$ for a fixed positive integer $k$. This would have implied the same statement of Theorem~\ref{theorem 3} with the sequence $ns$ replaced by $n^ks$. Unfortunately the argument fails for any $k\ge 2$ as then $H$ is no longer cocompact in $T'$ (due to the fact that the topology on $T'$ is not the product topology but the restricted one). Nonetheless, using a version of Theorem~\ref{unique ergodicity}
for the group $\SL_2$ (which is not available in the literature) and the choice $b_p=\diag{p,p^{-1}}$, leads to a proof of the validity of statement of Theorem~\ref{theorem 3} for the sequence $n^2s$. It seems
plausible  
that a better understanding of the proof of Theorem~\ref{unique ergodicity} could lead to a proof of the validity of the statement for $n^ks$ for general $k$ as well.    
\end{remark}
\begin{remark}\label{stronger form of theorem 3}
It is worth noting that a slight variant of the above argument actually yields a stronger uniform version
of Theorem~\ref{theorem 3} namely
\begin{theorem}\label{theorem 3'}
For any $M>0$ there exists a number $N$ such that for any irrational $s\in[0,1]$, there exists some $1\le n\le N$ for which $c(np)\ge M$. 
\end{theorem} 
\end{remark}
We end this section with two natural questions which emerge from the proof of Theorem~\ref{theorem 3}.
We use the notation presented in that proof. 
In the argument yielding the proof of Theorem~\ref{theorem 3} we used the assumption that the sequence 
$c(ns)$ is bounded to guarantee that the sequence of measures $\mu_n$ constructed in~\eqref{mu n} has no
escape of mass. It seems plausible that if the number $s$ is assumed to be badly approximable, then the non-escape of mass might be automatic for certain constructions of $\mu_n$. More precisely:
\begin{question}\label{natural question 1}
Is it true that for any badly approximable number $s\in[0,1]$, one can choose the F\o lner sequence $F_n'$ of $C'$, such that if $F_n$ is defined as in~\eqref{F_n}, with $m_n$ arbitrarily large, then the sequence 
of probability measures $\mu_n$ defined in~\eqref{mu n} has no escape of mass.
\end{question} 
We note that by applying the results from~\cite{AS} one can give a positive answer to Question~\ref{natural question 1} for quadratic irrationals which are of course badly approximable. It is not hard to see by applying Theorem~\ref{unique ergodicity}, that a positive answer for Question~\ref{natural question 1} leads to a positive answer to the following question:
\begin{question}\label{natural question 2}
Is it true that for any badly approximable number $s\in[0,1]$, and for any finite pattern $\textbf{w}=(w_1,\dots,w_\ell)
$ of natural numbers, there exists $n\in\bN$ such that the continued fraction expansion of $ns$ contains the pattern $\textbf{w}$ infinitely many times.
\end{question} 
\appendix 
\section{Proofs of several lemmas}\label{appendix}
In this section we give proofs for some of the lemmas appearing in \S\ref{preliminaries}.
\begin{proof}[Proof of Lemma ~\ref{WA vectors and dynamics}]
We think of points in $\PGL_{d+1}(\bR)/\PGL_{d+1}(\bZ)$ as unimodular lattices in $\bR^{d+1}$ as in \S\S~\ref{space of lattices}. For $v\in\bR^d$, the general form of a vector $w$ in the lattice $a_tu_v\PGL_{d+1}(\bZ)$ is given by
\begin{equation}\label{general vector}
w=\sum_1^d e^t(nv_i+m_i)+e^{-dt}ne_{d+1},
\end{equation}
where $e_i$ denotes the standard basis of $\bR^{d+1}$, $v_i$ denotes the $i$th' coordinate of $v$ and $m_i,n\in\bZ$.
Assume that $\eps>0$ is given so that the inequality 
\begin{equation}\label{bali}
\norm{nv+\vec{m}}_\infty <\frac{\eps}{n^{1/d}}
\end{equation}
has only finitely many solutions $\vec{m}\in\bZ^d,n\in\bZ\setminus\set{0}.$ We will show that for $w\ne 0$ as in ~\eqref{general vector} $\norm{w}_\infty>\eps$ for large enough $t$'s. Theorem ~\ref{Mahler} then implies the validity of the lemma. 

Let $N_0$ be given so that for $\av{n}\ge N_0$, there are no solutions to ~\eqref{bali}. For each $n$ with $0<\av{n}<N_0$, set $\del_{v,n}=\min_{\vec{m}\in\bZ^d}\norm{nv+\vec{m}}_\infty$, and for $n=0$ set
$\del_{v,0}=1.$ Note that as $v$ is irrational (otherwise there would have been infinitely many solutions
to ~\eqref{bali}), we have for all $ 0\le \av{n} < N_0$ that $\del_{v,n}>0$. We denote $\min_{0\le\av{n}<N_0} 
\del_{v,n}=\del$.  Let $T>0$ be such that for $t>T$, $e^t\del> 1$. Let $t>T$ be given. We now estimate the norm of 
$w\ne 0$ in ~\eqref{general vector}. There are two possibilities. If $0\le\av{n}<N_0$ then by construction, one of the 
first $d$ coordinates of $w$ is greater in absolute value than $e^t\del >1$. If $\av{n}\ge N_0$ then by the choice of $N_0$, ~\eqref{bali} is violated and there exists $1\le i\le d$ with $\av{nv_i+m_i}>\frac{\eps}{n^{1/d}}$. This means
that the product of the $i$-th coordinate of $w$ to the power of $d$, times the $(d+1)$-th coordinate satisfies
$\av{\pa{e^t(nv_i+m_i)}^{d}\pa{e^{-dt}n}}>\eps$. This shows (assuming $\eps<1$) that one of the coordinates of $w$ must be of absolute value greater than $\eps$, as desired. 
\end{proof} 
\begin{proof}[Proof of Lemma ~\ref{c.f.e}]
In this proof we use some basic facts about continued fractions. The reader is referred to
~\cite{EW} and to ~\cite{Vdp}. Let $s\in [0,1)$ be irrational with c.f.e.\ $s=[a_1,a_2\dots].$ For $n\in\bN$, let $p_n(s)=p_n,q_n(s)=q_n\in\bN$ be the co-prime positive
integers defined by the equation $p_n/q_n=[a_1,\dots a_n]$ (see ~\eqref{q.f}). $p_n/q_n$ is called the $n$-th convergent of $s$. The following two identities are well known
for all $n>0$: 
\begin{align}\label{identity}
q_{n+1}&=a_{n+1}q_n+q_{n-1},\\
s-\frac{p_n}{q_n}&=\sum_{k\ge n} (-1)^k\frac{1}{q_kq_{k+1}}.\nonumber
\end{align} 
It follows that $q_n\nearrow\infty$ and hence the above series is a Leibniz series and therefore we have  
\begin{equation}\label{i}
\av{s-\frac{p_n}{q_n}}\ge \frac{1}{q_nq_{n+1}}-\frac{1}{q_{n+1}q_{n+2}}=\frac{q_{n+2}-q_n}{q_nq_{n+1}q_{n+2}}=\frac{a_{n+2}}{q_nq_{n+2}},
\end{equation}
where the last equality follows from ~\eqref{identity}. By applying ~\eqref{identity} twice, we have
$q_{n+2}<(a_{n+2}+1)(a_{n+1}+1)q_n.$ This together with ~\eqref{i} yields
\begin{equation}\label{b}
\av{q_ns-p_n}\ge \frac{a_{n+2}}{(a_{n+2}+1)(a_{n+1}+1)q_n}.
\end{equation}
It is also well known that the convergents give the best possible approximations to $s$ in the following sense: For any rational $\frac{a}{b}$ with $0<b\le q_n$ one has $\av{q_ns-p_n}\le\av{bs-a}$. 
It follows that if
$c(s)=\limsup a_n$ satisfies $c(s)<N$ for some $N\in\bN$, then there are only finitely many solutions
$a,b\in\bZ, b\ne 0$, to the inequality
$$\av{bs+a}<\frac{(N+2)^{-2}}{b}.$$
Lemma ~\ref{WA vectors and dynamics} now gives us the desired result.
\end{proof}
For the proof of Lemma ~\ref{all patterns} we need some theory which we now survey. This theory dates back to the work of E. Artin (see ~\cite{Series}). For a thorough discussion we refer the reader to 
~\cite{EW}. We first note that $\PGL_2(\bR)/\PGL_2(\bZ)\simeq \PSL_2(\bR)/\PSL_2(\bZ)$. So we might as well carry on our analysis in the latter space. We let $\bH$ denote the upper half plane. 
On $\bH$ we take the Riemmannian metric defined as usual by taking at the tangent space to the point $z=x+iy\in\bH$, the inner product given by the usual Euclidean one, multiplied by $\frac{1}{y^2}$. With this
metric, the right action of $G=\PSL_2(\bR)$ on $\bH$ given by
\begin{equation}\label{action}
z\cdot\pa{
\begin{array}{ll}
a&b\\
c&d
\end{array}
}=\frac{dz-b}{-cz+a},
\end{equation}  
becomes an action by isometries. Hence, this action induces an action on the unit tangent bundle $T^1(\bH)$. One can easily check that this action is transitive and free, hence, once we choose a base point of  $T^1(\bH)$, the orbit map gives a diffeomorphism between $G$ and $T^1(\bH)$.
We make the common choice for the base point and choose the point $i^\uparrow$ which denotes the unit vector pointing upwards in the tangent space to $i\in\bH.$ Fixing this identification of $T^1(\bH)$ and $G$ once and for all, we are able to talk about the geodesic flow on $G$. 
It is an easy exercise to show that the geodesic flow is given by the action from the left of the diagonal group in $G$. More precisely, given $g\in G$ the point $a_{-t/2}g$ corresponds to the time $t$ flow starting at $g$. Hence the action of the group $a_t$ is then the \textit{backwards geodesic flow in double speed}. We define for each $g\in G$ the starting (resp. end) point of the geodesic through $g$, $e_-(g)$ (resp. $e_+(g)$), to be the intersection of the path $\set{a_tg}_{t>0}$, projected to $\bH$, (resp. $\set{a_tg}_{t<0}$) with the boundary of $\bH$ in $\bC\cup\set{\infty}$, namely with $\bR\cup\set{\infty}$.
In other words, in the notation of ~\eqref{action} we have 
\begin{equation}\label{starting}
e_-\pa{
\begin{array}{ll}
a&b\\
c&d
\end{array}
}=\frac{-b}{a},\quad
e_+\pa{
\begin{array}{ll}
a&b\\
c&d
\end{array}
}=\frac{-d}{c}
.
\end{equation}
We see that the starting point of $u_s$ is $s$. We now wish to connect the continued fraction expansion (c.f.e.)\ of $s$ with the geodesic ray $\set{a_tu_s}_{t>0}$ which starts at $s$.   
We denote the projection from $G$ to $G/\PSL_2(\bZ)$ by $\pi$. 
We will need the following three subsets of $G$:
\begin{align*}
C^+&=\set{g\in G:g \textrm{ lies on the } y \textrm{ axis, and } e_-(g)\in[0,1], e_+(g)<-1},\\ 
C^-&=\set{g\in G:g \textrm{ lies on the } y \textrm{ axis, and } e_-(g)\in[-1,0], e_+(g)>1},\\
C&=C^+\cup C^-. 
\end{align*}
The reader could prove the following theorem by simple geometric arguments (see ~\cite{EW}).
\begin{theorem}\label{Gauss}
The submanifold $C\subset G$ has the following properties:
\begin{enumerate}
\item $\pi:C\to \pi(C)$ is injective. Hence we have a canonical way of defining the starting (resp. end) point $e_-(x)$ (resp. $e_+(x)$) of $x\in \pi(C)$.
\item $\pi(C)$ is a cross section for the geodesic flow. We denote the first return  map (with respect to the $a_t$-action) by $\rho:\pi(C)\to\pi(C)$. A point $x\in \pi(C)$
returns to $\pi(C)$ infinitely often (i.e.\ $\rho^n(x)$ is defined for all $n>0$) if and only if $e_-(x)$ is irrational. In this case, its visits to $\pi(C)$ alternate between $\pi(C^+)$ and $\pi(C^-)$.
\item The map $x\mapsto \av{e_-(x)}$ from $\pi(C)$ to $[0,1]$ is a factor map connecting the first return map $\rho$
and the Gauss map on the unit interval (which is the shift on the c.f.e.).
\end{enumerate}  
\end{theorem}
The last bit of information we need in order to argue the proof
of Lemma ~\ref{all patterns}, is that if $s_1,s_2\in [0,1]\setminus \bQ$ satisfy $s_1=s_2\ga$ for some $\ga\in \PSL_2(\bZ)$ (the action given in ~\eqref{action}), then the continued fraction expansions of $s_1$ and $s_2$ only differ at their beginnings.
\begin{proof}[Proof of Lemma ~\ref{all patterns}]
Let $s\in[0,1]$ be such that $\set{a_tu_s\PSL_2(\bZ)}_{t>0}$ is dense in $G/\PSL_2(\bZ)$. In particular, $s$ is irrational. 
Given a pattern $(b_1,\dots ,b_k)\in\bN^k$, the set 
$$P=\set{s\in[0,1)\setminus\bQ:a_i(s)=b_i \textrm{ for } 1\le i\le k}$$ is an open set in $[0,1)\setminus \bQ$. It follows from (3) of Theorem ~\ref{Gauss}, that there is an open set $\tilde{P}\subset\pi(C)$, so that for any point $x\in\tilde{P}$, the starting point $e_-(x)$, if irrational, is in $P$. 
The density assumption gives us that there exists a sequence of times $t_i\nearrow\infty$ such that $x_i=a_{t_i}u_s\PGL_2(\bZ)\in\tilde{P}$ and 
moreover by (2) of Theorem ~\ref{Gauss} we may assume that $x_i\in C^+$, hence $e_-(x_i)\in[0,1]$. Now the c.f.e.\ of $s=e_-(u_s)$ differs from that of $e_-(x_1)$ only in their beginnings
(by the paragraph preceding this proof)
but by (3) of Theorem ~\ref{Gauss}, the c.f.e.\ of $e_-(x_1)$ must contain the pattern $b_1\dots b_k$ infinitely many times (as the c.f.e.\ of $e_-(x_i)$ starts with this
pattern and is a shift of the c.f.e.\ of $e_-(x_1)$), hence so does the c.f.e.\ of $s$ as desired.
\end{proof}
\def\cprime{$'$}

\end{document}